\documentclass[11pt,reqno]{amsart}
\usepackage{amsmath}
\usepackage{amssymb}
\usepackage{amstext}
\usepackage{amsfonts}
\usepackage{amsthm}
\usepackage{ifthen}
\usepackage{xspace}
\usepackage{comment}
\usepackage{graphicx}
\usepackage{mathrsfs}
\usepackage{enumerate}
\usepackage{breakcites}

\numberwithin{equation}{section}  

\newtheorem{punkt}{}[section]

\theoremstyle{plain}
\newtheorem{corollary}[punkt]{Corollary}
\newtheorem{lemma}[punkt]{Lemma}
\newtheorem{proposition}[punkt]{Proposition}
\newtheorem{theorem}[punkt]{Theorem}

\theoremstyle{definition}

\theoremstyle{plain}
\newtheorem*{corollary*}{Corollary}
\newtheorem*{lemma*}{Lemma}
\newtheorem*{proposition*}{Proposition}
\newtheorem*{theorem*}{Theorem}

\theoremstyle{definition}
\newtheorem*{remark*}{Remark}
\newtheorem*{remarks*}{Remarks}
\newtheorem*{example*}{Example}
\newtheorem*{examples*}{Examples}
\newtheorem*{definition*}{Definition}
\newtheorem*{conjecture*}{Conjecture}
\newtheorem*{assumption*}{Assumption}
\newtheorem*{assumptions*}{Assumptions}
\newtheorem*{construction*}{Construction}

\def\myo{\mathcal{O}}

\def\ee{\mathbb{E}}

\def\re{\qopname\relax{no}{Re}\,}
\def\im{\qopname\relax{no}{Im}\,}

\def\ai{\qopname\relax{no}{Ai}\nolimits}

\def\eg{e.g.\@\xspace}

\def\iid{i.i.d.\@\xspace}

\def\i{\qopname\relax{no}{\bf i}}

\def\eskip{\mskip24mu}

\begin{document}

\vspace*{4.0 \baselineskip}

\title{Asymptotics of Characteristic Polynomials \\[+5pt] of Wigner Matrices at the Edge of the Spectrum}

\author{H. K\"osters}
\address{Holger K\"osters, Fakult\"at f\"ur Mathematik, Universit\"at Bielefeld,
Postfach 100131, 33501 Bielefeld, Germany}
\email{hkoesters@math.uni-bielefeld.de}


\begin{abstract}
We investigate the asymptotic behaviour of the second-order correlation function
of the characteristic polynomial of a Hermitian Wigner matrix
at the edge of the spectrum. 
We show that the suitably rescaled second-order correlation function
is asymptotically given by the Airy kernel, thereby generalizing
the well-known result for the Gaussian Unitary Ensemble (GUE).
Moreover, we~obtain similar results for real-symmetric Wigner matrices.
\end{abstract}

\maketitle

\markboth{H. K\"osters}{Characteristic Polynomials of Wigner Matrices}
  
\section{Introduction and Statement of the Main Results}

Let $Q$ be a fixed probability distribution on the real line such that
\begin{align}
\label{moments-hermitian}
\int x \ Q(dx) = 0 \,,
\quad
a := \int x^2 \ Q(dx) = 1/2 \,,
\quad
b := \int x^4 \ Q(dx) < \infty \,,
\end{align}
and for any $N = 1,2,3,\hdots,$ 
let $X_N := (X_{ij})_{i,j=1,\hdots,N}$ denote 
the associated \linebreak Hermitian Wigner matrix of size $N$.
This means that
$$
X_{ij}
:= 
\begin{cases}
X^{\re}_{ij} + \i X^{\im}_{ij} & \text{ for } i < j , \\
\sqrt{2} X^{\re}_{ii}          & \text{ for } i = j , \\
X^{\re}_{ji} - \i X^{\im}_{ji} & \text{ for } i > j , \\
\end{cases}
$$
where $\{ X^{\re}_{ij} \,|\, i \leq j \} \cup \{ X^{\im}_{ij} \,|\, i < j \}$
is a collection of \iid real random variables with distribution $Q$.
The second-order correlation function of the characteristic poly\-nomial 
of the random matrix $X_N$ is defined by
$$
f_N(\mu,\nu) := \ee \big( D_N(\mu) \, D_N(\nu) \big) \qquad (\mu,\nu \in \mathbb{R}) \,,
$$
where $D_N(\lambda) := \det (X_N - \lambda)$.
We are interested in the asymptotic behaviour of $f_N(\mu_N,\nu_N)$
as $N \to \infty$, for certain sequences $(\mu_N)$, $(\nu_N)$
which will be specified below.
Furthermore, the correlation coefficient 
of the characteristic poly\-nomial of the random matrix $X_N$
is defined by 
$$
\sigma_N(\mu,\nu) := \frac{\ee \big( \big( D_N(\mu) - \ee D_N(\mu) \big) \, \big( D_N(\nu) - \ee D_N(\nu) \big) \big)}{\sqrt{ \ee \big( D_N(\mu) - \ee D_N(\mu) \big)^2 } \, \sqrt{ \ee \big( D_N(\nu) - \ee D_N(\nu) \big)^2 }} \qquad (\mu,\nu \in \mathbb{R}) \,.
$$

\pagebreak[2]
    
In the special case where $Q$ is the Gaussian distribution
with mean~$0$ and variance~$\frac{1}{2}$,
the distribution of the random matrix $X_N$
is the so-called Gaussian Unitary Ensemble (GUE)
(see \eg \textsc{Forrester} \cite{Fo} or \textsc{Mehta} \cite{Me}
but note that we work with a different variance).
In this case, it is well-known that 
\begin{align}
\label{link-hermitian}
f_N(\mu,\nu) = \sqrt{2\pi} \, N! \, e^{(\mu^2+\nu^2)/4} \, K_{N+1}(\mu,\nu) \,,
\end{align}
where 
$$
K_{N}(x,y) := e^{-(x^2+y^2)/4} \sum_{k=1}^{N} \frac{p_{k-1}(x) p_{k-1}(y)}{\sqrt{2\pi} (k-1)!}
$$
and the $p_k$ are the monic orthogonal polynomials
with respect to the weight function $e^{-x^2/2}$
(see \eg Chapter~4.1 in \textsc{Forrester} \cite{Fo}).
Thus, up to scaling, the $p_k$ coincide with the Hermite polynomials
(as defined in \textsc{Szeg\"o} \cite{Sz}),
and it is possible to derive the asymptotics 
of the second-order correlation function $f_N$ 
from the well-known asymptotics of the~Hermite polynomials
(see \eg Theorem~8.22.9 in \textsc{Szeg\"o} \cite{Sz}).
More precisely, one obtains the following (well-known) results
(see~also Chapter~4.2 in \textsc{Forrester} \cite{Fo}):
For $\xi \in (-2,+2)$ and any $\mu,\nu \in \mathbb{R}$,
\begin{align}
\label{bulk-GUE}
  \lim_{N \to \infty} c_N' \, f_N \left( \sqrt{N}\xi+\mu/\sqrt{N}\varrho(\xi), \sqrt{N}\xi+\mu/\sqrt{N}\varrho(\xi) \right)
= \mathbb{S}(\mu,\nu) \,,
\end{align}
where $c_N' := (\sqrt{2\pi} \, N! \, N^{1/2} \, \varrho(\xi) \, \exp(\frac{1}{2}N\xi^2 + \frac{1}{2}(\mu+\nu)\xi/\varrho(\xi)))^{-1}$,
$\varrho(\xi) := \tfrac{1}{2\pi} \sqrt{4 - \xi^2}$,
and
\begin{align}
\label{sine-kernel}
\mathbb{S}(\mu,\nu) := \frac{\sin \pi(\mu-\nu)}{\pi(\mu-\nu)} \,.
\end{align}
For $\xi = +2$ and any $\mu,\nu \in \mathbb{R}$,
\begin{align}
\label{edge-GUE}
  \lim_{N \to \infty} c_N'' \, f_N \left( 2\sqrt{N}+\mu/N^{1/6}, 2\sqrt{N}+\nu/N^{1/6} \right)
= \mathbb{A}(\mu,\nu) \,,
\end{align}
where $c_N'' :=(\sqrt{2\pi} \, N! \, N^{1/6} \, \exp(2N + (\mu+\nu)N^{1/3}))^{-1}$,
\begin{align}
\label{airy-kernel}
\mathbb{A}(\mu,\nu) := \frac{\ai(\mu) \ai'(\nu) - \ai'(\mu) \ai(\nu)}{\mu-\nu} \,,
\end{align}
and $\ai$ denotes the Airy function
(see \eg \textsc{Abramowitz} and \textsc{Stegun} \cite{AS}).
By~symmetry, a similar result holds for~$\xi = -2$.
The functions in (\ref{sine-kernel})~and~(\ref{airy-kernel})
are also called the sine kernel and the Airy kernel, respectively.
Furthermore, \linebreak it is well-known that the eigenvalues of a random matrix $X_N$ 
from the GUE are distributed roughly over the interval $[-2\sqrt{N},+2\sqrt{N}]$.
That~is~why the results (\ref{bulk-GUE}) and (\ref{edge-GUE})
are also said to refer to the bulk and the edge of the~spectrum, 
respectively.

Recently, \textsc{G\"otze} and \textsc{K\"osters} \cite{GK} have shown
that the result (\ref{bulk-GUE}) for the~bulk is (almost) ``universal''
in the sense that it holds not only for the GUE, 
but also (with minor modifications) for more general Hermitian Wigner matrices 
as introduced at the beginning of this section.
More precisely, under the assumption (\ref{moments-hermitian}),
we have
\begin{align}
\label{bulk-hermitian}
  \lim_{N \to \infty} c_N' \, f_N \left( \sqrt{N}\xi+\mu/\sqrt{N}\varrho(\xi), \sqrt{N}\xi+\mu/\sqrt{N}\varrho(\xi) \right)
= \exp(b\!-\!\tfrac{3}{4}) \, \mathbb{S}(\mu,\nu)
\end{align}
for all $\xi \in (-2,+2)$ and all $\mu,\nu \in \mathbb{R}$,
where $c_N'$, $\varrho(\xi)$ and $\mathbb{S}(\mu,\nu)$
are the same as in (\ref{bulk-GUE}).

Therefore, it seems natural to ask 
whether the result (\ref{edge-GUE}) for the edge 
can also be generalized to more general Hermitian Wigner matrices.
The main purpose of this paper is to answer this question
in the affirmative. More precisely, our first result is as follows:

\begin{theorem}
\label{theorem-hermitian}
Under (\ref{moments-hermitian}), we have
\begin{align}
\label{edge-hermitian}
  \lim_{N \to \infty} c_N'' \, f_N \left( 2\sqrt{N}+\mu/N^{1/6}, 2\sqrt{N}+\nu/N^{1/6} \right)
= \exp(b\!-\!\tfrac{3}{4}) \, \mathbb{A}(\mu,\nu)
\end{align}
for all $\mu,\nu \in \mathbb{R}$, where $c_N''$ and $\mathbb{A}(\mu,\nu)$ 
are the same as in (\ref{edge-GUE}).
\end{theorem}

\begin{corollary}
\label{corollary-hermitian}
Under (\ref{moments-hermitian}), we have
\begin{align}
\label{edge-hermitian2}
  \lim_{N \to \infty} \sigma_N \left( 2\sqrt{N}+\mu/N^{1/6}, 2\sqrt{N}+\nu/N^{1/6} \right)
= \frac{\mathbb{A}(\mu,\nu)}{\sqrt{\mathbb{A}(\mu,\mu)} \, \sqrt{\mathbb{A}(\nu,\nu)}}
\end{align}
for all $\mu,\nu \in \mathbb{R}$.
\end{corollary}  

Moreover, it turns out that similar results 
hold for real-symmetric Wigner matrices.
Let $\widetilde{Q}$ be a fixed probability distribution 
on the real line such that
\begin{align}
\label{moments-symmetric}
\int x \ \widetilde{Q}(dx) = 0 \,,
\quad
\widetilde{a} := \int x^2 \ \widetilde{Q}(dx) = 1 \,,
\quad
\widetilde{b} := \int x^4 \ \widetilde{Q}(dx) < \infty \,,
\end{align}
and for any $N = 1,2,3,\hdots,$ 
let $\widetilde{X}_N := (\widetilde{X}_{ij})_{i,j=1,\hdots,N}$ denote 
the associated \linebreak real-symmetric Wigner matrix of size $N$.
This means that
$$
\widetilde{X}_{ij}
:= 
\begin{cases}
\widetilde{X}^{\re}_{ij}          & \text{ for } i < j , \\
\sqrt{2} \widetilde{X}^{\re}_{ii} & \text{ for } i = j , \\
\widetilde{X}^{\re}_{ji}          & \text{ for } i > j , \\
\end{cases}
$$
where $\{ \widetilde{X}^{\re}_{ij} \,|\, i \leq j \}$ is a collection 
of \iid real random variables with distribution $\widetilde{Q}$.
Then, similarly as above, the second-order correlation function 
and the cor\-relation coefficient of the characteristic poly\-nomial 
of the random matrix $\widetilde{X}_N$ are defined by
$$
\widetilde{f}_N(\mu,\nu) := \ee \big( \widetilde{D}_N(\mu) \, \widetilde{D}_N(\nu) \big) \qquad (\mu,\nu \in \mathbb{R})
$$
and
$$
\widetilde\sigma_N(\mu,\nu) := \frac{\ee \big( \big( \widetilde{D}_N(\mu) - \ee \widetilde{D}_N(\mu) \big) \, \big( \widetilde{D}_N(\nu) - \ee \widetilde{D}_N(\nu) \big) \big)}{\sqrt{ \ee \big( \widetilde{D}_N(\mu) - \ee \widetilde{D}_N(\mu) \big)^2 } \, \sqrt{ \ee \big( \widetilde{D}_N(\nu) - \ee \widetilde{D}_N(\nu) \big)^2 }} \qquad (\mu,\nu \in \mathbb{R}) \,,
$$
respectively, where now $\widetilde{D}_N(\lambda) := \det(\widetilde{X}_N - \lambda)$.

Following the approach by \textsc{G\"otze} and \textsc{K\"osters} \cite{GK},
\textsc{K\"osters} \cite{Ko} recently showed that  
under the assumption (\ref{moments-symmetric}),
we have
\begin{align}
\label{bulk-symmetric}
  \lim_{N \to \infty} d_N' \, \widetilde{f}_N \left( \sqrt{N}\xi+\mu/\sqrt{N}\varrho(\xi), \sqrt{N}\xi+\mu/\sqrt{N}\varrho(\xi) \right)
= \exp(\tfrac{\widetilde{b}-3}{2}) \, \mathbb{T}(\mu,\nu)
\end{align}
for all $\xi \in (-2,+2)$ and all $\mu,\nu \in \mathbb{R}$, \pagebreak[1]
where $\varrho(\xi)$ is the same as in (\ref{bulk-GUE}), \linebreak
$d_N' := (\sqrt{2\pi} \, N! \, N^{3/2} \, \varrho(\xi)^3 \, \exp(\frac{1}{2}N\xi^2 + \frac{1}{2}(\mu+\nu)\xi/\varrho(\xi)))^{-1}$,
and
\begin{align}
\label{sine-kernel-2}
\mathbb{T}(\mu,\nu) := \frac{2 \sin \pi(\mu-\nu)}{\pi (\mu-\nu)^3} - \frac{2 \cos \pi(\mu-\nu)}{(\mu-\nu)^2} \,.
\end{align}
Thus, we also have universality (in the same sense as above)
in the bulk of real-symmetric Wigner matrices.
In the special case where $\widetilde{Q}$ 
is the Gaussian distribution with mean~$0$ and variance~$1$,
the distribution of the random matrix $\widetilde{X}_N$
is the so-called Gaussian Orthogonal Ensemble (GOE)
(see \eg \textsc{Forrester} \cite{Fo} or \textsc{Mehta} \cite{Me}
but, again, note that we work with a different variance).
In this case, (\ref{bulk-symmetric}) had been obtained previously
by \textsc{Br\'ezin} and \textsc{Hikami} \cite{BH2}.

Our second result shows that the result (\ref{bulk-symmetric}) 
admits an analogue for the edge of the spectrum, too:
\begin{theorem}
\label{theorem-symmetric}
Under (\ref{moments-symmetric}), we have
\begin{align}
\label{edge-symmetric}
  \lim_{N \to \infty} d_N'' \, \widetilde{f}_N \left( 2\sqrt{N}+\mu/N^{1/6}, 2\sqrt{N}+\nu/N^{1/6} \right)
= \exp(\tfrac{\widetilde{b}-3}{2}) \, \mathbb{B}(\mu,\nu)
\end{align}
for all $\mu,\nu \in \mathbb{R}$, where
$d_N'' :=(\sqrt{2\pi} \, N! \, N^{1/2} \, \exp(2N + (\mu+\nu)N^{1/3}))^{-1}$,
and
\begin{align}
\label{airy-kernel-2}
\mathbb{B}(\mu,\nu) := \frac{(\mu+\nu) \ai(\mu) \ai(\nu) - 2 \ai'(\mu) \ai'(\nu)}{(\mu-\nu)^2} + \frac{2 \ai(\mu) \ai'(\nu) - 2 \ai'(\mu) \ai(\nu)}{(\mu-\nu)^3} \,.
\end{align}
\end{theorem}

\begin{corollary}
\label{corollary-symmetric}
Under (\ref{moments-symmetric}), we have
\begin{align}
\label{edge-symmetric2}
  \lim_{N \to \infty} \widetilde\sigma_N \left( 2\sqrt{N}+\mu/N^{1/6}, 2\sqrt{N}+\nu/N^{1/6} \right)
= \frac{\mathbb{B}(\mu,\nu)}{\sqrt{\mathbb{B}(\mu,\mu)} \, \sqrt{\mathbb{B}(\nu,\nu)}}
\end{align}
for all $\mu,\nu \in \mathbb{R}$.
\end{corollary}  

\noindent{}In the special case of the GOE, (\ref{edge-symmetric}) 
can already be found in \textsc{Br\'ezin}~and~\textsc{Hikami} \cite{BH3}.

\pagebreak[2]

It seems interesting to note that the functions $\mathbb{S}$ and $\mathbb{T}$
arising for the bulk of the spectrum are related by the identity
$$
\mathbb{T}(x,y) = \left( \frac{1}{x-y} \left( \frac{\partial}{\partial y} - \frac{\partial}{\partial x} \right) \right) \mathbb{S}(x,y)
$$
and that the functions $\mathbb{A}$ and $\mathbb{B}$
arising for the edge of the spectrum are related by the analogous identity
$$
\mathbb{B}(x,y) = \left( \frac{1}{x-y} \left( \frac{\partial}{\partial y} - \frac{\partial}{\partial x} \right) \right) \mathbb{A}(x,y) \,,
$$
see also \textsc{Br\'ezin} and \textsc{Hikami} \cite{BH3}.
(To check the latter identity, use the fact that the Airy function $\ai(z)$ 
satisfies the differential equation $\ai''(z) = z \ai(z)$.)

Also, observe that in all the cases previously mentioned, the precise choice
of the underlying distribution $Q$ or $\widetilde{Q}$
enters into the asymptotic behaviour of the second-order
correlation function of the characteristic polynomial only
as a multiplicative factor depending on the fourth cumulant.
Thus, the results are essentially ``universal''.

Let us mention some related results from the literature.
It is well-known that at the edge of the spectrum of Wigner matrices,
we have universality also for the cor\-relation function of the eigenvalues 
themselves (see \textsc{Soshnikov} \cite{So}).
In~contrast to~that, for the bulk of the spectrum of Wigner matrices,
only partial results are available in this direction
(see \textsc{Johansson} \cite{Jo}).
Furthermore, in the special~cases of the GUE and the GOE,
several authors have investigated the averages
of more general products and even ratios of characteristic polynomials
(see \eg \textsc{Br\'ezin} and \textsc{Hikami} \cite{BH1,BH2,BH3},
\textsc{Fyodorov} and \textsc{Strahov} \cite{FS3}, 
\textsc{Baik}, \textsc{Deift} and \textsc{Strahov} \cite{BDS}, 
\textsc{Strahov} and \textsc{Fyodorov} \cite{SF3}, 
\textsc{Akemann} and \textsc{Fyodorov} \cite{AF3},
\textsc{Vanlessen} \cite{Va},
\textsc{Borodin} and \textsc{Strahov} \cite{BS}).
Even more, at least in the~Hermitian setting,
some of these results have been shown to be ``universal''
in that they continue to hold (with some modifications) 
for the class of unitary-invariant ensembles.
%
%
For Wigner matrices, however, less seems to be known
in this~respect.

\medskip

\textbf{Acknowledgement:} I would like to thank Alexander Soshnikov 
for the suggestion to study the problem at the edge of the spectrum.

\bigskip

\section{Outline of the Proofs}

To prove Theorems \ref{theorem-hermitian} and \ref{theorem-symmetric},
we will start from the fact that for fixed $\mu,\nu \in \mathbb{R}$,
the exponential generating functions of the sequences 
$f_N(\mu,\nu)$ and $\widetilde{f}_N(\mu,\nu)$
are given explicitly by
$$
\sum_{N=0}^{\infty} f_N(\mu,\nu) \frac{x^N}{N!}
=
\frac{\exp \left( \mu \nu \cdot \frac{x}{1-x^2} - \tfrac{1}{2} (\mu^2 + \nu^2) \cdot \frac{x^2}{1-x^2} + \big( b\!-\!\frac{3}{4} \big) x^2 \right)}{(1-x)^{3/2} \cdot (1+x)^{1/2}} \qquad (|x| < 1)
$$
(see Lemma~2.3 in \textsc{G\"otze} and \textsc{K\"osters} \cite{GK})
and
$$
\sum_{N=0}^{\infty} \widetilde{f}_N(\mu,\nu) \frac{x^N}{N!}
=
\frac{\exp \left( \mu \nu \cdot \frac{x}{1-x^2} - \tfrac{1}{2} (\mu^2 + \nu^2) \cdot \frac{x^2}{1-x^2} + \big( \frac{\widetilde{b}-3}{2} \big) x^2 \right)}{(1-x)^{5/2} \cdot (1+x)^{1/2}} \qquad (|x| < 1)
$$
(see Lemma 2.3 in \textsc{K\"osters} \cite{Ko}), respectively.
This fact opens up the possibility to study 
the asymptotic behaviour of the second-order correlation functions
by evaluating appropriate contour integrals 
of their exponential generating functions.
In fact, this strategy was already used 
by \textsc{G\"otze} and \textsc{K\"osters} \cite{GK} and \textsc{K\"osters} \cite{Ko}
to~obtain the above-mentioned results for the bulk of the spectrum.
Here we carry out a similar analysis for the edge of the spectrum.

Since it does not require any additional efforts,
it seems convenient to evaluate the values
\begin{align}
\label{fa-definition}
\frac{f_N^{(\alpha)}(\mu,\nu)}{N!} := \frac{1}{2\pi i} \int_\gamma \frac{\exp \left( \mu \nu \cdot \frac{z}{1-z^2} - \tfrac{1}{2} (\mu^2 + \nu^2) \cdot \frac{z^2}{1-z^2} + b^* z^2 \right)}{(1-z)^{\alpha + (1/2)} \cdot (1+z)^{1/2}} \ \frac{dz}{z^{N+1}}
\end{align}
for arbitrary $\alpha > 0$ and $b^* \in \mathbb{R}$,
where $\gamma$ denotes a contour around the origin.
By the~fore\-going, the case $\alpha = 1$ corresponds to the Hermitian case
and the case $\alpha = 2$ corresponds to the real-symmetric case.

We remark in passing that for a general parameter $\alpha > 0$,
it is not hard to see that the exponential generating function under consideration
can be interpreted as that of the second-order correlation~function
of the characteristic polynomial of a random matrix 
from (a rescaled version of) the tridiagonal beta ensemble
introduced by \textsc{Dumitriu} and \textsc{Edelman} \cite{DE1}
with $\alpha = 2/\beta$. For this~inter\-pretation, one should set $b^* := 0$.

For $\alpha > 0$ and $\mu,\nu \in \mathbb{R}$, put
\begin{align}
\label{idefinition}
I^{(\alpha)}(\mu,\nu) := \frac{1}{4 \pi^{3/2}} \int_{-\infty}^{+\infty} \frac{\exp \big( \tfrac{1}{12}(1\!-\!iu)^3 - \tfrac{1}{2}(\mu+\nu)(1\!-\!iu) - \tfrac{1}{4}(\mu-\nu)^2 / (1\!-\!iu) \big)}{(1-iu)^{\alpha+(1/2)}} \ du \,.
\end{align}
In the next section, we will prove the following results:

\begin{proposition}
\label{theproposition2}
For any $\alpha > 0$, $b^* \in \mathbb{R}$ and $\mu,\nu \in \mathbb{R}$, we have
\begin{multline*}
\lim_{N \to \infty} \left( \sqrt{2\pi} \, N! \, N^{(2\alpha-1)/6} \exp(2N+(\mu+\nu)N^{1/3}) \right)^{-1} \\ \,\cdot\, f_N^{(\alpha)}(2\sqrt{N}+\mu N^{-1/6},2\sqrt{N}+\nu N^{-1/6})
=
\exp(b^*) \, I^{(\alpha)}(\mu,\nu) \,.
\end{multline*}
\end{proposition}

\begin{proposition}
\label{theproposition3}
For any $\alpha \in \mathbb{N}$ and $\mu,\nu \in \mathbb{R}$, we have
\begin{align*}
I^{(\alpha)}(\mu,\nu)
=
\left( \frac{1}{\mu-\nu} \left( \frac{\partial}{\partial \nu} - \frac{\partial}{\partial \mu} \right) \right)^{(\alpha)} \Big( \ai(\mu) \ai(\nu) \Big) \,,
\end{align*}
where $(\,\cdot\,)^{(\alpha)}$ denotes the $\alpha$-fold application
of the given differential operator. \linebreak
(For $\mu = \nu$, consider the continuous extension of the right-hand side.)
\end{proposition}

It is straightforward to check that
$$
I^{(1)}(x,y) = \left( \frac{1}{x-y} \left( \frac{\partial}{\partial y} - \frac{\partial}{\partial x} \right) \right)^{(1)} \Big( \ai(x) \ai(y) \Big) = \mathbb{A}(x,y) \,\phantom{.}
$$
and \mbox{\qquad\qquad}
$$
I^{(2)}(x,y) = \left( \frac{1}{x-y} \left( \frac{\partial}{\partial y} - \frac{\partial}{\partial x} \right) \right)^{(2)} \Big( \ai(x) \ai(y) \Big) = \mathbb{B}(x,y) \,.
$$
(For the second identity, use the fact that the Airy function $\ai(z)$ 
satisfies the differential equation $\ai''(z) = z \ai(z)$.)
Thus, in view of the preceding comments, it is immediate that 
Theorems \ref{theorem-hermitian} and \ref{theorem-symmetric},
which correspond to the special cases $\alpha = 1$ and $\alpha = 2$,
follow from Propositions \ref{theproposition2} and \ref{theproposition3}.

To prove Corollaries \ref{corollary-hermitian} and \ref{corollary-symmetric},
observe that 
\begin{align*}
   \sigma_N(\mu,\nu) 
&= \frac{f_N(\mu,\nu) - \ee D_N(\mu) \, \ee D_N(\nu)}{\sqrt{ f_N(\mu,\mu) - \big( \ee D_N(\mu) \big)^2 } \, \sqrt{ f_N(\nu,\nu) - \big( \ee D_N(\nu) \big)^2 }}
\end{align*}
in the Hermitian case and
\begin{align*}
   \widetilde\sigma_N(\mu,\nu) 
&= \frac{\widetilde{f}_N(\mu,\nu) - \ee \widetilde{D}_N(\mu) \, \ee \widetilde{D}_N(\nu)}{\sqrt{ \widetilde{f}_N(\mu,\mu) - \big( \ee \widetilde{D}_N(\mu) \big)^2 } \, \sqrt{ \widetilde{f}_N(\nu,\nu) - \big( \ee \widetilde{D}_N(\nu) \big)^2 }}
\end{align*}
in the real-symmetric case. Moreover, it is not difficult to see that
$$
\ee D_N(\lambda) = \ee \widetilde{D}_N(\lambda) = g_N(\lambda) := (-1)^N \, 2^{-N/2} \, H_N(\lambda/\sqrt{2}) \qquad (\lambda \in \mathbb{R}) \,,
$$ 
where $H_N(x)$ is the $N$th Hermite polynomial \pagebreak[1]
(see \eg Section~5.5 in \textsc{Szeg\"o} \cite{Sz}).  
Thus, in both cases, the expectation of the characteristic polynomial
is given by the~\emph{same}~function $g_N(\lambda)$.
Therefore, to deduce Corollaries \ref{corollary-hermitian} and \ref{corollary-symmetric}
from Theorems \ref{theorem-hermitian} and \ref{theorem-symmetric}, respectively,
it will be sufficient to show that $g_N(\mu) \, g_N(\nu)$ is~asymptotically negligible
in comparison to $f_N(\mu,\nu)$ and $\widetilde{f}_N(\mu,\nu)$.

Slightly more generally, we will consider the case of an arbitrary parameter $\alpha > 0$
and investigate the asymptotic behaviour of
\begin{align}
\label{sigma-definition}
    \sigma^{(\alpha)}_N(\mu,\nu) 
&:= \frac{f^{(\alpha)}_N(\mu,\nu) - g_N(\mu) \, g_N(\nu)}{\sqrt{ f^{(\alpha)}_N(\mu,\mu) - g_N(\mu)^2 } \, \sqrt{ f^{(\alpha)}_N(\nu,\nu) - g_N(\nu)^2 }}
\end{align}
for any $\mu,\nu \in \mathbb{R}$. In the next section, 
we will show that Proposition \ref{theproposition2}
entails the~following result:

\begin{proposition}
\label{theproposition5}
For any $\alpha > 0$, $b^* \in \mathbb{R}$ and $\mu,\nu \in \mathbb{R}$
such that $I^{(\alpha)}(\mu,\mu) > 0$, $I^{(\alpha)}(\nu,\nu) > 0$, we have
$$
\lim_{N \to \infty} \sigma_N^{(\alpha)}(2\sqrt{N}+\mu N^{-1/6},2\sqrt{N}+\nu N^{-1/6})
\\ =
\frac{I^{(\alpha)}(\mu,\nu)}{\sqrt{I^{(\alpha)}(\mu,\mu)} \, \sqrt{I^{(\alpha)}(\nu,\nu)}} \,.
$$
\end{proposition}

Furthermore, we will prove the following:

\begin{proposition}
\label{theproposition6}
For any $\alpha \in \mathbb{N}$, we have $I^{(\alpha)}(x,x) > 0$
for all $x \in \mathbb{R}$.
\end{proposition}

In view of the preceding comments, it is obvious that 
Corollaries \ref{corollary-hermitian} and \ref{corollary-symmetric},
which correspond to the special cases $\alpha = 1$ and $\alpha = 2$, 
follow from Propositions \ref{theproposition5}~and~\ref{theproposition6}.

We remark in passing that for a general parameter $\alpha > 0$,
$\sigma^{(\alpha)}_N(\mu,\nu)$ can be interpreted as 
the correlation coefficient of the characteristic polynomial of a random matrix 
from (the rescaled version of) the tridiagonal beta ensemble
(with $\alpha = 2/\beta$), since in this setting, the average
of the characteristic polynomial is also given by the Hermite polynomial
(see Theorem~4.1 in \textsc{Dumitriu} and \textsc{Edelman} \cite{DE1}).

\bigskip

\section{The Proofs}

\begin{proof}[Proof of Proposition \ref{theproposition2}]
Fix $\mu,\nu \in \mathbb{R}$.
We have to evaluate
\begin{align}
\label{intrep1}
\frac{f^{(\alpha)}_N(\mu_N,\nu_N)}{N!}
=
\frac{1}{2 \pi i} \int_\gamma \frac{\exp \left( \mu_N \nu_N \cdot \frac{z}{1-z^2} - \tfrac{1}{2} (\mu_N^2 + \nu_N^2) \cdot \frac{z^2}{1-z^2} + b^* z^2 \right)}{(1-z)^{\alpha+(1/2)} \cdot (1+z)^{1/2}} \ \frac{dz}{z^{N+1}} \,,
\end{align}
where $\mu_N := 2N^{1/2} + \mu N^{-1/6}$, $\nu_N := 2N^{1/2} + \nu N^{-1/6}$,
and $\gamma$ denotes a contour around the origin
(which will be chosen further below).

Setting $\xi_N := (\mu_N+\nu_N)/2$ and $\eta_N := (\mu_N-\nu_N)/2$, 
(\ref{intrep1}) may be written as
\begin{equation}
\label{intrep2}
\exp \left( \tfrac{1}{2} \xi_N^2 + \eta_N^2 \right) \,\cdot\, \frac{1}{2 \pi i} \int_\gamma \frac{\exp \left( - \tfrac{1}{2} \xi_N^2 \cdot \frac{1-z}{1+z} - \eta_N^2 \cdot \frac{1}{1-z} + b^* z^2  \right)}{(1-z)^{\alpha+(1/2)} \cdot (1+z)^{1/2}} \ \frac{dz}{z^{N+1}}
\end{equation}
with 
$$
\xi_N^2 = 4N + 2(\mu+\nu) N^{1/3} + \tfrac{1}{4}(\mu+\nu)^2 N^{-1/3}
\qquad\text{and}\qquad
\eta_N^2 = \tfrac{1}{4} (\mu-\nu)^2 N^{-1/3} \,.
$$
In particular, the leading exponential factor in (\ref{intrep2}) 
is asymptotically the same as~that in Proposition \ref{theproposition2}.

Thus, to complete the proof of Proposition \ref{theproposition2},
it suffices to show that
\begin{multline}
\label{intrep3}
\lim_{N \to \infty} \frac{N^{-(2\alpha-1)/6}}{(2 \pi)^{3/2} i} \int_\gamma \frac{\exp \left( - \tfrac{1}{2} \xi_N^2 \cdot \frac{1-z}{1+z} - \eta_N^2 \cdot \frac{1}{1-z} + b^* z^2  \right)}{(1-z)^{\alpha+(1/2)} \cdot (1+z)^{1/2}} \ \frac{dz}{z^{N+1}}
\\ =
\frac{\exp(b^*)}{4 \pi^{3/2}} \int_{-\infty}^{+\infty} \frac{\exp(h_\infty(1-iu))}{(1-iu)^{\alpha+(1/2)}} \ du \,, \quad
\end{multline}
where
\begin{align}
\label{hinfty}
h_\infty(z) := \tfrac{1}{12}z^3 - \tfrac{1}{2}(\mu+\nu)z - \tfrac{1}{4}(\mu-\nu)^2 / z \,.
\end{align}

Similarly as in the proof of the main theorem in \cite{GK}, 
the basic idea is that the~main contribution to the integral 
in (\ref{intrep3}) comes from a~small neighborhood of the point $z=1$.
Let
$$
h_N(z) := - \tfrac{1}{2} \xi_N^2 \cdot \frac{1-z}{1+z} - \eta_N^2 \cdot \frac{1}{1-z} + b^* z^2 - (\alpha + \tfrac{1}{2}) \log (1-z) - \tfrac{1}{2} \log (1+z)  - N \log z
$$
(where $\log$ denotes the principal branch of the logarithm)
and
$$
\gamma_N(t) := (1 - N^{-1/3}) \exp(it) \,, \qquad -\pi \leq t \leq +\pi \,.
$$
Then the left-hand side in (\ref{intrep3}) can be rewritten as
\begin{multline*}
  \frac{N^{-(2\alpha-1)/6}}{(2 \pi)^{3/2} i} \int_{\gamma_N} \frac{\exp \left( - \tfrac{1}{2} \xi_N^2 \cdot \frac{1-z}{1+z} - \eta_N^2 \cdot \frac{1}{1-z} + b^* z^2  \right)}{(1-z)^{\alpha+(1/2)} \cdot (1+z)^{1/2}} \ \frac{dz}{z^{N+1}} \\
= \frac{N^{-(2\alpha-1)/6}}{(2 \pi)^{3/2} i} \int_{\gamma_N} \exp(h_N(z)) \ \frac{dz}{z}
= \frac{N^{-(2\alpha-1)/6}}{(2 \pi)^{3/2}} \int_{-\pi}^{+\pi} \exp(h_N(\gamma_N(t))) \ dt \,.
\end{multline*}
Put $I_{N,1}(a) := (-a N^{-1/3}; +a N^{-1/3})$
and $I_{N,2}(a) := (-\pi,+\pi) \setminus (-a N^{-1/3}; +a N^{-1/3})$,
where $a > 0$. We shall prove the following:

\medskip

\noindent\emph{Claim 1:}
For any fixed $a > 0$,
\begin{multline*}
\lim_{N \to \infty} \frac{N^{-(2\alpha-1)/6}}{(2\pi)^{3/2}} \int_{I_{N,1}(a)} \exp(h_N(\gamma_N(t))) \ dt
=
\frac{\exp(b^*)}{4\pi^{3/2}} \int_{-a}^{+a} \frac{\exp(h_\infty(1-iu))}{(1-iu)^{\alpha+(1/2)}} \ du \,. 
\end{multline*}

\medskip

\noindent\emph{Claim 2:}
For any $\delta > 0$, there exists a constant $a_0 > 0$ 
such that for all $a \geq a_0$,
$$
\left| \frac{N^{-(2\alpha-1)/6}}{(2\pi)^{3/2}} \int_{I_{N,2}(a)} \exp(h_N(\gamma_N(t))) \ dt \right| \leq \delta
$$
for all $N \in \mathbb{N}$ sufficiently large.

\medskip

Before turning to the proofs, let us show that Claims 1~and~2 imply (\ref{intrep3}).
Observe~that
$$
\int_{-\infty}^{+\infty} \frac{\exp(h_\infty(1-iu))}{(1-iu)^{\alpha+(1/2)}} \ du < \infty \,.
$$

\pagebreak[3]

\noindent{}Thus, for $a > 0$ sufficiently large, we have 
not only the conclusion of Claim~2,
but also the inequality
$$
\left| \frac{\exp(b^*)}{4\pi^{3/2}} \int_{\mathbb{R} \setminus (-a,+a)} \frac{\exp(h_\infty(1-iu))}{(1-iu)^{\alpha+(1/2)}} \ du \right| < \delta \,.
$$
Hence, in combination with Claim~1, it follows that
\begin{align*}
&\eskip \left| \frac{N^{-(2\alpha-1)/6}}{(2 \pi)^{3/2}} \int_{-\pi}^{+\pi} \exp(h_N(\gamma_N(t))) \ dt - \frac{\exp(b^*)}{4\pi^{3/2}} \int_{-\infty}^{+\infty} \frac{\exp (h_\infty(1-iu))}{(1-iu)^{\alpha+(1/2)}} \ du \right| \\
&\leq \, \left| \frac{N^{-(2\alpha-1)/6}}{(2 \pi)^{3/2}} \int_{I_{N,2}(a)} \exp(h_N(\gamma_N(t))) \ dt \right| \\
& \qquad\,+\, \left| \frac{N^{-(2\alpha-1)/6}}{(2\pi)^{3/2}} \int_{I_{N,1}(a)} \exp(h_N(\gamma_N(t))) \ dt - \frac{\exp(b^*)}{4\pi^{3/2}} \int_{-a}^{+a} \frac{\exp (h_\infty(1-iu))}{(1-iu)^{\alpha+(1/2)}} \ du \right| \\
& \qquad\,+\, \left| \frac{\exp(b^*)}{4\pi^{3/2}} \int_{\mathbb{R} \setminus (-a,+a)} \frac{\exp (h_\infty(1-iu))}{(1-iu)^{\alpha+(1/2)}} \ du \right| \\
&\leq
  \delta + \delta + \delta
= 3\delta
\end{align*}
for all $N \in \mathbb{N}$ sufficiently large.
Since $\delta > 0$ is arbitrary, this proves (\ref{intrep3}).

\medskip

\emph{Proof of Claim 1:}
First of all, substituting $t = u N^{-1/3}$, we have
$$
\frac{N^{-(2\alpha-1)/6}}{(2 \pi)^{3/2}} \int_{-a N^{-1/3}}^{+a N^{-1/3}} \exp(h_N(\gamma_N(t))) \ dt
= 
\frac{N^{-(2\alpha+1)/6}}{(2 \pi)^{3/2}} \int_{-a}^{+a} \exp(h_N(\gamma_N(u N^{-1/3}))) \ du \,.
$$
We will determine the asymptotics of $h_N(\gamma_N(u N^{-1/3}))$ as $N \to \infty$, 
for \mbox{$u \in [-a,+a]$}.
The $\myo$-bounds occurring in the sequel hold uniformly in this region.
To begin with, note that
\begin{align*}
   z := \gamma_N(u N^{-1/3}) 
&= (1-N^{-1/3}) e^{iu N^{-1/3}} \\ 
&= (1-N^{-1/3}) (1 + iu N^{-1/3} - \tfrac{1}{2} u^2 N^{-2/3} - \tfrac{1}{6} iu^3 N^{-1} + \myo(N^{-4/3})) \\
&= 1 - (1-iu) N^{-1/3} - (iu + \tfrac{1}{2} u^2) N^{-2/3} + (\tfrac{1}{2} u^2 - \tfrac{1}{6} iu^3) N^{-1} + \myo(N^{-4/3}) \,.
\end{align*}
Therefore, for $N$ sufficiently large, we have the following approximations:
\begin{align*}
&\eskip {-} \tfrac{1}{2} \xi_N^2 \cdot \frac{1-z}{1+z} \\
&= - \tfrac{1}{2} \xi_N^2 \Big( \tfrac{1}{2}(1-z) + \tfrac{1}{4}(1-z)^2 + \tfrac{1}{8}(1-z)^3 + \myo((1-z)^4) \Big) \\
&= - \tfrac{1}{2} \left( 4N + 2(\mu+\nu) N^{1/3} + \myo(N^{-1/3}) \right) \\&\qquad\,\cdot\, \left( \tfrac{1}{2} (1-iu) N^{-1/3} + \tfrac{1}{4} N^{-2/3} + (\tfrac{1}{8}+\tfrac{1}{8}iu+\tfrac{1}{8}u^2-\tfrac{1}{24}iu^3) N^{-1} + \myo(N^{-4/3}) \right) \\
&= - (1-iu) N^{2/3} - \tfrac{1}{2} N^{1/3} + \Big( \tfrac{1}{12}(1-iu)^3 - \tfrac{1}{2}(\mu+\nu)(1-iu) - \tfrac{1}{3} \Big) + \myo(N^{-1/3})
\end{align*}
\begin{align*}
   - \eta_N^2 \cdot \frac{1}{1-z}
&= - \left( \tfrac{1}{4} (\mu-\nu)^2 N^{-1/3} \right) \cdot \frac{1}{(1-iu) N^{-1/3}} \cdot \left( 1 + \myo(N^{-1/3}) \right) \\
&= - \tfrac{1}{4} (\mu-\nu)^2 / (1-iu) \cdot \left( 1 + \myo(N^{-1/3}) \right) 
\end{align*}
$$
b^* z^2 = b^* + \myo(N^{-1/3})
$$
$$
- (\alpha+\tfrac{1}{2}) \log (1-z) = + \tfrac{1}{3} (\alpha+\tfrac{1}{2}) \log N - (\alpha + \tfrac{1}{2}) \log (1 - i u) + \myo(N^{-1/3})
$$
$$
- \tfrac{1}{2} \log (1+z) = - \tfrac{1}{2} \log 2 + \myo(N^{-1/3})
$$
\begin{align*}
   - N \log z 
&= - N \left( iuN^{-1/3} + \log (1-N^{-1/3}) \right) \\
&= - N \left( iuN^{-1/3} - N^{-1/3} - \tfrac{1}{2} N^{-2/3} - \tfrac{1}{3} N^{-1} + \myo(N^{-4/3}) \right) \\
&= (1 - iu) N^{2/3} + \tfrac{1}{2} N^{1/3} + \tfrac{1}{3} + \myo(N^{-1/3})
\end{align*}
Putting these approximations together, the terms of highest order cancel out,
and we end up with
\begin{multline*}
h_N(\gamma_N(u N^{-1/3})) = \tfrac{1}{3} \left( \alpha+\tfrac{1}{2} \right) \log N + \Big( \tfrac{1}{12}(1-iu)^3 - \tfrac{1}{2}(\mu+\nu)(1-iu) \\ - \tfrac{1}{4}(\mu-\nu)^2 / (1-iu) - (\alpha+\tfrac{1}{2}) \log (1-iu) + b^* - \tfrac{1}{2} \log 2 \Big) + \myo(N^{-1/3}) \,.
\end{multline*}
Hence, by an application of the dominated convergence theorem,
it follows that
\begin{multline*}
\lim_{N \to \infty} \frac{N^{-(2\alpha+1)/6}}{(2\pi)^{3/2}} \int_{-a}^{+a} \exp(h_N(\gamma_N(u N^{-1/3}))) \ du \\
=
\frac{\exp(b^*)}{4\pi^{3/2}} \int_{-a}^{+a} \frac{\exp \Big( \tfrac{1}{12}(1-iu)^3 - \tfrac{1}{2}(\mu+\nu)(1-iu) - \tfrac{1}{4}(\mu-\nu)^2 / (1-iu) \Big)}{(1-iu)^{\alpha+(1/2)}} \ du \,,
\end{multline*}
and Claim~1 is proved.

\pagebreak[2]
\medskip

\emph{Proof of Claim 2:}
By symmetry, it suffices to consider the interval $(a N^{-1/3},\pi)$.
Write the integral as
$$
\frac{N^{-(2\alpha-1)/6}}{(2 \pi)^{3/2}} \int_{a N^{-1/3}}^{\pi} \frac{\exp \left( - \tfrac{1}{2} \xi_N^2 \cdot \frac{1-\gamma_N(t)}{1+\gamma_N(t)} - \eta_N^2 \cdot \frac{1}{1-\gamma_N(t)} + b^* \gamma_N(t)^2  \right)}{(1-\gamma_N(t))^{\alpha+(1/2)} \cdot (1+\gamma_N(t))^{1/2}} \ \frac{dt}{\gamma_N(t)^{N}} \,.
$$
Since $\xi_N$ and $\eta_N$ are real and $|\gamma_N(t)| \leq 1$, this is clearly bounded by
$$
\frac{N^{-(2\alpha-1)/6}}{(2 \pi)^{3/2}} \int_{a N^{-1/3}}^{\pi} \frac{\exp \left( - \tfrac{1}{2} \xi_N^2 \cdot \re\left(\frac{1-\gamma_N(t)}{1+\gamma_N(t)}\right) - \eta_N^2 \cdot \re\left(\frac{1}{1-\gamma_N(t)}\right) + |b^*| \right)}{|1-\gamma_N(t)|^{\alpha+(1/2)} \cdot |1+\gamma_N(t)|^{1/2}} \ \frac{dt}{|\gamma_N(t)|^{N}} \,.
$$
In the following, let $K,K_1,K_2 > 0$ denote constants
which depend only on $\alpha,b^*,\mu,\nu$
(which are regarded as fixed) but which may change 
from occurrence to occurrence. Then we have
$$
\xi_N^2 \geq 4N - K N^{1/3}
$$
and
$$
\re \left(\frac{1-\gamma_N(t)}{1+\gamma_N(t)}\right) = \frac{2N^{-1/3} - N^{-2/3}}{(2 + 2 \cos t)(1 - N^{-1/3}) + N^{-2/3}}
$$
(as follows from a straightforward calculation)
and therefore
$$
\tfrac{1}{2} \xi_N^2 \re \left(\frac{1-\gamma_N(t)}{1+\gamma_N(t)}\right) 
\geq 
\frac{4N^{2/3} - 2N^{1/3} - K}{(2 + 2 \cos t)(1 - N^{-1/3}) + N^{-2/3}} \,.
$$
Since
\begin{align*}
   |\gamma_N(t)|^N
&= \exp(N \log (1 - N^{-1/3})) \\
&\geq \exp(N (- N^{-1/3} - \tfrac{1}{2} N^{-2/3} - K N^{-1})) \\
&= \exp(- N^{2/3} - \tfrac{1}{2} N^{1/3} - K) \,,
\end{align*}
it follows that
\begin{align*}
&\eskip \frac{\exp \left( -\tfrac{1}{2} \xi_N^2 \re\left(\frac{1-\gamma_N(t)}{1+\gamma_N(t)}\right) \right)}{|\gamma(t)|^N} \\ 
&\leq \exp \left( - \frac{4N^{2/3} - 2N^{1/3} - K_1}{(2+2\cos t)(1-N^{-1/3}) + N^{-2/3}} + N^{2/3} + \tfrac{1}{2} N^{1/3} + K_2 \right) \\
&\leq \exp \left( - \frac{(2-2\cos t)N^{2/3} - (1-\cos t)N^{1/3} - K^*}{(2+2\cos t)(1-N^{-1/3}) + N^{-2/3}} \right) 
\end{align*}
for some constant $K^* > 0$, say. Now let $\varepsilon > 0$ denote a constant
such that $\varepsilon^2 t^2 \leq 1 - \cos t \leq \tfrac{1}{2} t^2$
for all $t \in [0,\pi]$.
Then, for $a \geq \varepsilon^{-1} \sqrt{2K^*}$,
$N^{1/3} \geq \max \{ \varepsilon^{-2},a\pi^{-1} \}$ and $t \in [a N^{-1/3},\pi]$, 
we~have
\begin{multline*}
     (2-2\cos t)N^{2/3} - (1-\cos t)N^{1/3} - K^* \\
\geq \varepsilon^2 t^2 N^{2/3} + ( \tfrac{1}{2} \varepsilon^2 t^2 N^{2/3} - \tfrac{1}{2} t^2 N^{1/3} ) + ( \tfrac{1}{2} \varepsilon^2 t^2 N^{2/3} - K^* )
\geq N^{2/3} \varepsilon^2 t^2 
\end{multline*}
and therefore
\begin{multline*}
     \frac{\exp \left( -\tfrac{1}{2} \xi_N^2 \re\left(\frac{1-\gamma_N(t)}{1+\gamma_N(t)}\right) \right)}{|\gamma(t)|^N} \\
\leq \exp \left( - \frac{N^{2/3} \varepsilon^2 t^2}{(2+2\cos t)(1-N^{-1/3}) + N^{-2/3}} \right)
\leq \exp \left( - N^{2/3} \varepsilon^2 t^2 / 4 \right) \,.
\end{multline*}
Thus, since $|1 \pm \gamma_N(t)| \geq N^{-1/3}$,
the integral under consideration is bounded by
\begin{align*}
&\eskip K \, N^{-(2\alpha-1)/6} \, \int_{a N^{-1/3}}^{\pi} \frac{\exp(- N^{2/3} \varepsilon^2 t^2 / 4)}{|1-\gamma_N(t)|^{\alpha+1/2} \cdot |1+\gamma_N(t)|^{1/2}} \ dt \\ 
&\leq
K \, N^{-(2\alpha-1)/6} \, \int_{a N^{-1/3}}^{\pi} \frac{\exp(- N^{2/3} \varepsilon^2 t^2 / 4)}{N^{-(2\alpha+1)/6}} \ dt \\
&=
K \, N^{1/3} \, \int_{a N^{-1/3}}^{\pi} \exp(- N^{2/3} \varepsilon^2 t^2 / 4) \ dt \\ 
&\leq
K \, \int_{a}^{\infty} \exp(- \varepsilon^2 u^2 / 4) \ du \,.
\end{align*}
Obviously, this upper bound can be made arbitrarily small 
by picking $a$ and $N$ large enough. This proves Claim~2.

\medskip

The proof of Proposition \ref{theproposition2} is complete now.
\end{proof}

\pagebreak[2]

The basic ingredient for the proof of Proposition \ref{theproposition3}
will be the following integral representation for the product 
of two Airy functions:

\begin{lemma}
\label{airyproduct}
For any $x,y \in \mathbb{C}$,
\begin{align}
\ai(x) \ai(y) = \frac{1}{4 \pi^{3/2} i} \int_{\mathcal{L}} \frac{\exp \big( \tfrac{1}{12}z^3 - \tfrac{1}{2}(x+y)z - \tfrac{1}{4}(x-y)^2 / z \big)}{z^{1/2}} \ dz \,,
\end{align}
where $\mathcal{L}$ denotes some (unbounded) contour
from $\infty e^{-\pi i/3}$ to $\infty e^{+\pi i/3}$.
\end{lemma}

In the special case where $y = \pm x$, 
this result can already be found in \textsc{Reid} \cite{Re}.
For the convenience of the reader, we give a detailed proof
of Lemma \ref{airyproduct}:

\begin{proof}
We start from the following well-known integral representation of the Airy function:
$$
\ai(x) = \frac{1}{2 \pi i} \int_{\mathcal{L}} \exp(\tfrac{1}{3} z^3 - xz) \ dz \,.
$$
A standard application of Cauchy's theorem shows that
the contour $\mathcal{L}$ can be deformed into the contour
$t \mapsto 1 + it$, $t \in \mathbb{R}$. Thus, we obtain
\begin{align}
\label{ai-int}
\ai(x) = \frac{1}{2 \pi i} \int_{1-i\infty}^{1+i\infty} \exp(\tfrac{1}{3} z^3 - xz) \ dz \,.
\end{align}
Observe that the resulting integral exists in the Lebesgue sense,
since we have
$$
\left| \exp( \tfrac{1}{3}(1+it)^3 - x(1+it) ) \right|
\leq 
\exp ( \tfrac{1}{3} - t^2 + |x|(1+|t|) )
$$
for any $t \in \mathbb{R}$. 

\pagebreak[1]

It follows from (\ref{ai-int}) that
$$
\ai(x) \ai(y) = \frac{1}{4 \pi^2} \iint \exp(\tfrac{1}{3} (1+it)^3 - x(1+it)) \exp(\tfrac{1}{3} (1+iu)^3 - y(1+iu)) \ du \ dt \,.
$$
Substituting $(t,u) = (\tfrac{1}{2}(v+w),\tfrac{1}{2}(v-w))$
and doing a small calculation, we find that
$$
\ai(x) \ai(y) = \frac{1}{8 \pi^2} \iint \exp( \tfrac{1}{12}(2+iv)^3 - \tfrac{1}{4}(2+iv)w^2 - \tfrac{1}{2}(x+y)(2+iv) - \tfrac{1}{2}(x-y)iw ) \ dw \ dv \,.
$$
Using the well-known relation
$$
\int_{-\infty}^{+\infty} \exp(-aw^2 -bw) \ dw = \frac{\sqrt{\pi}}{\sqrt{a}} \, e^{b^2/4a} 
$$
(where $\re(a) > 0$), it follows that
$$
\ai(x) \ai(y)
=
\frac{1}{4\pi^{3/2}} \int \frac{\exp( \tfrac{1}{12}(2+iv)^3 - \tfrac{1}{2}(x+y)(2+iv) - \tfrac{1}{4}(x-y)^2 / (2+iv))}{(2+iv)^{1/2}} \ dv
$$
or
$$
\ai(x) \ai(y)
=
\frac{1}{4\pi^{3/2}i} \int_{2-i\infty}^{2+i\infty} \frac{\exp( \tfrac{1}{12}z^3 - \tfrac{1}{2}(x+y)z - \tfrac{1}{4}(x-y)^2 / z)}{z^{1/2}} \ dv \,.
$$
By another application of Cauchy's theorem,
the contour $v \mapsto 2 + iv$, $v \in \mathbb{R}$,
may be deformed back into the contour $\mathcal{L}$.
\end{proof}

\begin{proof}[Proof of Proposition \ref{theproposition3}]
Replacing the contour $\mathcal{L}$ in Lemma \ref{airyproduct} 
by the contour $t \mapsto 1 + it$ and~substituting $t = -u$, 
we have
$$
\ai(x) \ai(y) = \frac{1}{4 \pi^{3/2}} \int_{-\infty}^{+\infty} \frac{\exp \big( \tfrac{1}{12}(1\!-\!iu)^3 - \tfrac{1}{2}(x+y)(1\!-\!iu) - \tfrac{1}{4}(x-y)^2 / (1\!-\!iu) \big)}{(1-iu)^{1/2}} \ du \,.
$$
By means of abbreviation, write $E(x,y,u)$ for the numerator inside the integral.
Then, for any $\alpha > 0$, we have
$$
\frac{\partial}{\partial y} \left( \int \frac{E(x,y,u)}{(1-iu)^{\alpha}} \ du \right)
=
\int \frac{E(x,y,u)}{(1-iu)^{\alpha}} \left( -\tfrac{1}{2}(1-iu) + \tfrac{1}{2}(x-y)/(1-iu) \right) \ du \,,
$$
$$
\frac{\partial}{\partial x} \left( \int \frac{E(x,y,u)}{(1-iu)^{\alpha}} \ du \right)
=
\int \frac{E(x,y,u)}{(1-iu)^{\alpha}} \left( -\tfrac{1}{2}(1-iu) - \tfrac{1}{2}(x-y)/(1-iu) \right) \ du \,,
$$
and therefore
$$
\left( \frac{1}{x-y} \left( \frac{\partial}{\partial y} - \frac{\partial}{\partial x} \right) \right) \left( \int \frac{E(x,y,u)}{(1-iu)^{\alpha}} \ du \right)
=
\int \frac{E(x,y,u)}{(1-iu)^{\alpha+1}} \ du \,.
$$
The assertion of Proposition \ref{theproposition3} now follows by induction.
\end{proof}

%
%
%
%
%

\begin{proof}[Proof of Proposition \ref{theproposition5}] \ 
Fix $\mu,\nu \in \mathbb{R}$, and put 
$\mu_N := 2N^{1/2} + \mu N^{-1/6}$, \linebreak $\nu_N := 2N^{1/2} + \nu N^{-1/6}$.
Using well-known results about the asymptotic properties of the Hermite polynomials
(see \eg Theorem 8.22.9\,(c) in \textsc{Szeg\"o} \cite{Sz}), we find that the function 
$g_N(\lambda)$ given by $g_N(\lambda) := (-1)^N \, 2^{-N/2} \, H_N(\lambda/\sqrt{2})$
($\lambda \in \mathbb{R}$) satisfies
\begin{align*}
&\eskip g_N(2\sqrt{N}+\mu N^{-1/6}) \, g_N(2\sqrt{N}+\nu N^{-1/6}) \\ 
&= 2^{-N} \, H_N(\sqrt{2N}+\mu N^{-1/6} / \sqrt{2}) \, H_N(\sqrt{2N}+\nu N^{-1/6} / \sqrt{2}) \\
&= \myo \left( 2^{-N} \, \exp(N + \mu N^{1/3}) \, 2^{N/2} \, N!^{1/2} \, N^{-1/12} \, \exp(N + \nu N^{1/3}) \, 2^{N/2} \, N!^{1/2} \, N^{-1/12} \right) \\
&= \myo \left( N! \, N^{-1/6} \, \exp(2N+(\mu+\nu)N^{1/3}) \right) \\
&= o \left( N! \, N^{(2\alpha-1)/6} \, \exp(2N+(\mu+\nu)N^{1/3}) \right)
\end{align*}
for any $\alpha > 0$. 
Setting $c_N^{(\alpha)}(\mu,\nu) := \left( \sqrt{2\pi} \, N! \, N^{(2\alpha-1)/6} \exp(2N+(\mu+\nu)N^{1/3}) \right)^{-1}$,
we therefore obtain
\begin{align*}
&\eskip \lim_{N \to \infty} \sigma_N^{(\alpha)}(2\sqrt{N}+\mu N^{-1/6},2\sqrt{N}+\nu N^{-1/6}) \\
&= \lim_{N \to \infty} \frac{f^{(\alpha)}_N(\mu_N,\nu_N) - g_N(\mu_N) \, g_N(\nu_N)}{\sqrt{ f^{(\alpha)}_N(\mu_N,\mu_N) - g_N(\mu_N)^2 } \, \sqrt{ f^{(\alpha)}_N(\nu_N,\nu_N) - g_N(\nu_N)^2 }} \\
&= \lim_{N \to \infty} \frac{c_N^{(\alpha)}(\mu,\nu) \left( f^{(\alpha)}_N(\mu_N,\nu_N) - g_N(\mu_N) \, g_N(\nu_N) \right)}{\sqrt{ c_N^{(\alpha)}(\mu,\mu) \left( f^{(\alpha)}_N(\mu_N,\mu_N) - g_N(\mu_N)^2 \right) } \, \sqrt{ c_N^{(\alpha)}(\nu,\nu) \left( f^{(\alpha)}_N(\nu_N,\nu_N) - g_N(\nu_N)^2 \right) }} \\
&= \frac{\lim_{N \to \infty} c_N^{(\alpha)}(\mu,\nu) \, f^{(\alpha)}_N(\mu_N,\nu_N)}{\sqrt{ \lim_{N \to \infty} c_N^{(\alpha)}(\mu,\mu) \, f^{(\alpha)}_N(\mu_N,\mu_N) } \, \sqrt{ \lim_{N \to \infty} c_N^{(\alpha)}(\nu,\nu) \, f^{(\alpha)}_N(\nu_N,\nu_N) }} \\
&= \frac{I^{(\alpha)}(\mu,\nu)}{\sqrt{I^{(\alpha)}(\mu,\mu)} \, \sqrt{I^{(\alpha)}(\nu,\nu)}} \,,
\end{align*}
where we have used Proposition \ref{theproposition2} as well as 
the assumptions $I^{(\alpha)}(\mu,\mu) > 0$, $I^{(\alpha)}(\nu,\nu)> 0$.
This~completes the proof of Proposition \ref{theproposition5}.
\end{proof}

\begin{proof}[Proof of Proposition \ref{theproposition6}]
First of all, note that the definition (\ref{idefinition}) 
may be extended to the case $\alpha = 0$ 
and that $I^{(0)}(x,x) = \ai(x)^2$ for any $x \in \mathbb{R}$
by~Lemma~\ref{airyproduct}. Thus, $I^{(0)}(x,x) \geq 0$ 
for any $x \in \mathbb{R}$, with strict inequality for $x > 0$
(since it is well-known that the Airy function does not have 
any zeroes on the positive half-axis).~More\-over, 
note that for any $\alpha > 0$,
\begin{align*}
   \frac{\partial}{\partial x} \bigg( I^{(\alpha)}(x,x) \bigg)
&= \frac{\partial}{\partial x} \bigg( \frac{1}{4 \pi^{3/2}} \int_{-\infty}^{+\infty} \frac{\exp \big( \tfrac{1}{12}(1\!-\!iu)^3 - x(1\!-\!iu) \big)}{(1-iu)^{\alpha+(1/2)}} \ du \bigg) \\
&= - \frac{1}{4 \pi^{3/2}} \int_{-\infty}^{+\infty} \frac{\exp \big( \tfrac{1}{12}(1\!-\!iu)^3 - x(1\!-\!iu) \big)}{(1-iu)^{\alpha-(1/2)}} \ du \\
&= - I^{(\alpha-1)}(x,x)
\end{align*}
for any $x \in \mathbb{R}$. Since $\lim_{x \to \infty} I^{(\alpha)}(x,x) = 0$,
this implies that for any $\alpha > 0$,
\begin{align*}
I^{(\alpha)}(x,x) = \int_{x}^{\infty} I^{(\alpha-1)}(y,y) \ dy
\end{align*}
for any $x \in \mathbb{R}$. Proposition \ref{theproposition6} now follows 
by a straightforward induction on~$\alpha$.
\end{proof}

\newpage

\bigskip

\bigskip

\end{document}